\documentclass[a4paper,11pt,reqno]{amsart}
\usepackage[a4paper,hscale=0.7,vscale=0.75,centering]{geometry}
\usepackage{mathrsfs} 

\usepackage{amsthm}
\usepackage{hyperref}
\usepackage{amssymb}
\usepackage{latexsym}
\usepackage{amsmath}
\usepackage{amsfonts}
\usepackage[dvips]{graphicx}
\usepackage[utf8]{inputenc}

\usepackage[LGR,T1]{fontenc}%
\usepackage[english]{babel}
\usepackage{url}
\usepackage{enumerate}
\usepackage{hyperref}
\usepackage{color}
\usepackage{bbm}
\usepackage{fancyhdr}
\newtheorem{Definition}{Definition}[part]
\newtheorem{Proposition}{Proposition}[part]

\newtheorem{Lemma}{Lemma}[part]
\newtheorem{Corollary}{Corollary}[part]
\newtheorem{Remark}{Remark}[part]

\makeatletter \@addtoreset{equation}{section}

\@addtoreset{Definition}{section}

\@addtoreset{Theorem}{section}

\@addtoreset{Proposition}{section}

\@addtoreset{Property}{section}

\@addtoreset{Assumption}{section}

\@addtoreset{Corollary}{section}

\@addtoreset{Lemma}{section}

\@addtoreset{Remark}{section}

\@addtoreset{Example}{section}

\newcommand{\cA}{\mathcal{A}}

\newcommand{\cF}{\mathcal{F}}

\newcommand{\cY}{\mathcal{Y}}
\newcommand{\cZ}{\mathcal{Z}}

\newcommand{\E}{\mathbb{E}}
\newcommand{\F}{\mathbb{F}}

\renewcommand{\P}{\mathbb{P}}

\newcommand{\R}{\mathbb{R}}

\newcommand{\wh}[1]{\ensuremath{ \widehat{#1} }}

\def \proof{{\noindent \bf Proof. }}

\newcommand{\ud}{\mathrm{d}}

\newcommand{\HYP}[1]
    {\ensuremath{({H#1} ) }}

\newcommand{\set}[1]
    {\ensuremath{\{ #1 \}}}
    
\newcommand{\HP}[1] %L2DPDT sur 0,T
    {\ensuremath{\mathscr{H}^{#1}}}

\newcommand{\SP}[1]{\ensuremath{\mathscr{S}^{#1}}}

\newcommand{\esp}[1]{\ensuremath{\mathbb{E} \!\! \left[#1\right] }}

\newcommand{\EFp}[2]
    {\ensuremath{% raccourci esperance conditionnelle discretisation
     \mathbb{E}_{#1}\!\!\left[#2\right] }}
     
\newcommand{\ti}[1]{t_{i #1}}

\renewcommand{\Xi}[1]{X_{i #1}}
\newcommand{\Yi}[1]{Y_{i #1}}
\newcommand{\Zi}[1]{Z_{i #1}}

\newcommand{\whEFp}[2]
    {\ensuremath{%
     \widehat{\mathbb{E}}_{#1}\!\left[#2\right] }}

\newcommand{\abs}[1]{\left\vert#1\right\vert}

\def\sqw{\hbox{\rlap{\leavevmode\raise.3ex\hbox{$\sqcap$}}$%
\sqcup$}}
\def\sqb{\hbox{\hskip5pt\vrule width4pt height6pt depth1.5pt%
\hskip1pt}}

\def\qed{\ifmmode\hbox{\hfill\sqb}\else{\ifhmode\unskip\fi%
\nobreak\hfil
\penalty50\hskip1em\null\nobreak\hfil\sqb
\parfillskip=0pt\finalhyphendemerits=0\endgraf}\fi}
\def\cqfd{\ifmmode\sqw\else{\ifhmode\unskip\fi\nobreak\hfil
\penalty50\hskip1em\null\nobreak\hfil\sqw
\parfillskip=0pt\finalhyphendemerits=0\endgraf}\fi}
\title{ Numerical stability analysis of the Euler scheme for BSDEs}

\author{Jean-Fran\c{c}ois Chassagneux}
\address[Jean-Fran\c{c}ois Chassagneux]{Department of Mathematics, Imperial College London, 180 Queen's Gate, London, SW7 2AZ, United Kingdom.}
\email{j.chassagneux@imperial.ac.uk}
\thanks{The work of the first author was partially supported by the Research Grant ANR-11-JS01-0007 - LIQUIRISK and EPSRC Mathematics Platform grant EP/I019111/1. } 
\author{Adrien Richou}
\address[Adrien Richou]{Univ. Bordeaux, IMB, UMR 5251, F-33400 Talence, France.}
\email{adrien.richou@math.u-bordeaux1.fr}

\begin{document}
\bibliographystyle{plain}

\begin{abstract}
In this paper, we study  the qualitative behaviour of approximation schemes for Backward Stochastic Differential Equations (BSDEs) by introducing a new notion of numerical stability.
For the Euler scheme, we provide sufficient conditions in the one-dimensional and multidimensional case to guarantee the numerical stability. We then perform a classical Von Neumann stability analysis in the case of a linear driver $f$ and exhibit necessary conditions to get stability in this case. Finally, we illustrate our results with numerical applications.
\end{abstract}

\maketitle

\vspace{5mm}


\vspace{5mm}

\noindent {\bf MSC Classification (2000):}  93E20, 65C99, 60H30.

\tableofcontents

 \newpage

%\input{stabBSDE}
%!TEX root = main.tex

\section{Introduction}
In this paper, we study  the qualitative behaviour of a class of numerical methods for Backward Stochastic Differential Equations (BSDEs) by introducing a new notion of numerical stability. Even though we will focus exclusively on the numerical schemes, we recall,
to motivate our work the definition of BSDEs in a classical setting, see e.g. \cite{Pardoux-Peng-90}. Let $(\Omega,\cA,\P)$ be a probability space supporting a $d$-dimensional Brownian motion $(W_t)_{t \ge 0}$. We denote by $\F:=(\cF_t)_{t\ge 0}$ the Brownian filtration.
Let $T>0$, $\xi$ be an $\cF_T$-measurable and square-integrable random variable and  $f: \Omega \times \R^+\times \R^d \times \R^{m \times d}\rightarrow \R^m$ in such a way that the process $(f(t,y,z))_{t \in [0,T]}$ is progressively measurable for all $(y,z) \in \R^m\times\R^{m\times d}$ and $\E \left[ \int_0^T \abs{f(s,0,0)}^2 ds\right]<+\infty$. The solution $(\cY,\cZ)$ of a BSDEs satisfies
\begin{align} \label{eq BSDE}
\cY_t = \xi + \int_t^T f(s,\cY_s, \cZ_s)\ud s - \int_t^T \cZ_s \ud W_s\;.
\end{align}
%\textcolor{red}{Mettre nouvelle notation car on utilise peu et diffencie avec schema}

If we assume that $f$ is a Lipschitz function with respect to $y$ and $z$ then it is known \cite{Pardoux-Peng-90} that the BSDE \eqref{eq BSDE} has a unique solution $(\cY,\cZ) \in \SP{2} \times \HP{2}$
where $\SP{2}$ is the set of continuous adapted processes satisfying $\esp{\sup_{s\in [0,T]}|U_s|^2}<\infty$ and $\HP{2}$ is the set
of progressively measurable processes $V$ satisfying $\esp{\int_0^T |V_t|^2 \ud t}<\infty$. Let us mention also that it is possible to relax some assumptions on $f$ and $\xi$: see e.g. \cite{Pardoux-99} for monotone generators with respect to $y$, \cite{Briand-Delyon-Hu-Pardoux-Stoica-03} for $L^p$ solutions and \cite{Kobylanski-00} for quadratic generators with respect to $z$.
These equations have applications e.g. in PDE analysis through non-linear Feynman-Kac formula \cite{Pardoux-Peng-92, cridel12}, stochastic control theory \cite{majyon99} or mathematical finance \cite{ElKaroui-Peng-Quenez-97}.
Recently, they have been used as non linear pricing methods \cite{cre12a,cre12b,bripav14,briliu14}.
In the past ten years, a lot of work has also been done on the numerical approximation of the above BSDE (and extensions) see e.g. \cite{Zhang-04,Bouchard-Touzi-04,Gobet-Labart-07,Chassagneux-13} and the references therein, especially in a markovian setting. This means that the terminal condition and the random part of the generator are given by deterministic measurable functions of a forward diffusion $X$, precisely $\xi := g(X_T)$ and $f(t,y,z)=\bar{f}(t,X_t,y,z)$, with
\begin{align*}
 X_t = X_0 + \int_0^tb(X_s)\ud s + \int_0^t \sigma(X_s) \ud W_s\;,\; t\in[0,T]\,.
\end{align*}
Here, we assume that $g$, $x \mapsto \bar f(t,x,y,z)$, $b$ and $\sigma$ are Lipschitz-continuous function.

One of the first numerical method that has been proposed, see e.g.  \cite{Zhang-04,Bouchard-Touzi-04} and the references therein for early works, is given by a discrete backward programming equation. Given a grid $\pi=\set{t_0:=0,\dots,t_i,\dots,t_n := T}$, one sets $Y_n = \xi$ and compute at each step:
\begin{align}
 \Yi{} &= \EFp{\ti{}}{\Yi{+1}+ (\ti{+1}-\ti{})f(t_i,\Yi{},\Zi{})} \label{eq th scheme Y imp}
 \\
 \Zi{} &= \EFp{\ti{}}{\frac1{\ti{+1}-\ti{}}\Yi{+1}(\Delta W_{i})'} \label{eq th scheme Z}
\end{align}
where $\Delta W_i := W_{\ti{+1}}-W_{\ti{}}$, $\EFp{t}{\cdot}$ stands for $\esp{\cdot \,| \, \cF_t}$ and $'$ denotes the transpose operator.

The above scheme is implicit in $Y$ and one can compute alternatively 
\begin{align}
\Yi{} = \EFp{\ti{}}{\Yi{+1}+ (\ti{+1}-\ti{})f(t_{i},\Yi{+1},\Zi{})} \label{eq th scheme Y exp}
\end{align}
to get an explicit version.

It has been shown that, in the Lipschitz setting, the above method has order at least one-half \cite{Zhang-04,Bouchard-Touzi-04} and generally at most one \cite{Gobet-Labart-07,Chassagneux-Crisan-13}.
Recently, other methods have been proposed with high order of convergence, one step methods as Runge-Kutta method \cite{Chassagneux-Crisan-13} and linear multi-step method \cite{Chassagneux-13}. A first motivation for this work comes from the need to distinguish between 'good' and 'bad' methods provided by the above papers. Indeed, the order of convergence of a scheme is an asymptotic property that allows to classify schemes when the number of time-steps tends to infinity. We would like here to know the quality of a scheme when the number of timesteps is set.

The study we perform in this paper can be seen as an extension of the numerical stability study of numerical schemes for ODEs, in the context of BSDEs. Indeed, if one considers a deterministic terminal value for $\xi$ and a deterministic generator $f$, in the one-dimensional setting the BSDE reduces to an ODE, and the corresponding scheme to an Euler method for ODE. It is well known, see e.g. \cite{Butcher-03}, that implicit and explicit Euler method have different stability behaviour in practice when $f$ is monotone. In particular, the explicit Euler method may become unstable if the timestep $h$ is not small enough. One should expect this behaviour also for the above BSDE scheme. The framework of BSDEs is (in some sense) richer than the one for ODEs and studying the stability of the methods given in \eqref{eq th scheme Y imp}-\eqref{eq th scheme Z} or \eqref{eq th scheme Y exp}-\eqref{eq th scheme Z} is already a challenging task. Moreover, to the best of our knowledge, it is the first time that such a 
study is undertaken. In the next paragraph, we motivate our work with an example belonging 'purely' to the BSDEs framework.

\subsection{A motivating example}
\label{subse motivation}
Let us consider  the BSDE \eqref{eq BSDE} with the following choice of coefficients: $g(\cdot)=\cos(\alpha \cdot)$, $X=W$ (dimension one) and $f(t,y,z)=bz$, for given real numbers $\alpha$ and $b$. Namely, $(\cY,\cZ)$ is solution to, 
\begin{align}
\cY_t = cos(\alpha W_T) + \int_t^T b\cZ_s \ud s - \int_t^T \cZ_s \ud W_s\;, t \le T\;.
\end{align}
At time $t=0$, the $\cY$ component is easily computed and given by
\begin{align*}
\cY_0 = e^{-\alpha^2\frac{T}{2}}\cos(\alpha b T)\;.
\end{align*}

We observe that $\cY_0$ is bounded by $1$, the bound of the terminal condition and moreover, $\cY_0 \rightarrow 0$ as $T \rightarrow \infty$.
This can be interpreted as a stability property of the BSDE, see next section, Proposition \ref{pr estim BSDE} .

We then consider the numerical approximation introduced in \eqref{eq th scheme Y imp}-\eqref{eq th scheme Z} above. In order to compute the conditional expectations and set the terminal value, we simply use a trinomial (recombining) tree for the Brownian motion, see e.g. \cite{Chassagneux-13} and an equidistant time grid of $[0,T]$ with $h=\frac{T}{n}$ and $n+1$ time steps. It is well known that, in this context, the error for the $\cY$ part is given by
%\begin{align*}
$
\cY_0 - Y_0^{} = O(h),
$
%\end{align*}
where $Y_0$ is the solution at time $0$ returned by the scheme, see \cite{Gobet-Labart-07,Chassagneux-13}. %This means that the scheme is convergent with an order equal to $1$.\\
Let us now try to observe this behaviour in practice by plotting the error $|\cY_0 - {Y}_0^{}|$ against the number of step, in logarithmic scale, for different value of $(b,T)$, $\alpha$ being set to $1$.

\begin{figure}[htb!] 
\includegraphics[height = 6cm, width = 
\textwidth]{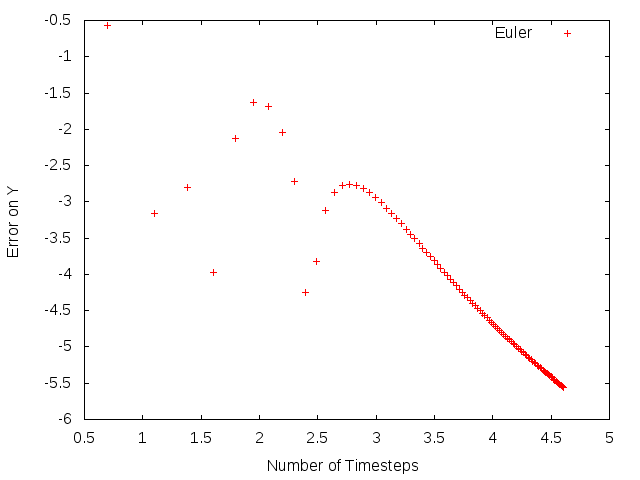} %[width = \textwidth]
\caption{Euler Scheme, $(b=1,T=10)$} \label{fig intro cz1T10N100}
\end{figure}

\begin{figure}[htb!] 
\includegraphics[height = 6cm, width = 
\textwidth]{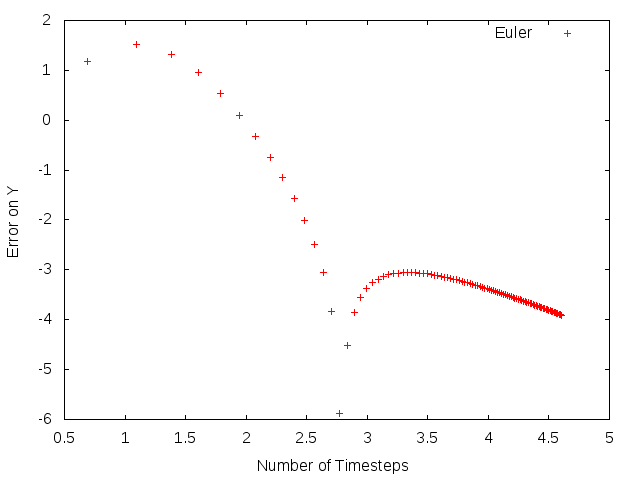} %[width = \textwidth]
\caption{Euler Scheme, $(b=5,T=1)$} \label{fig intro cz5T1N100}
\end{figure}

\begin{figure}[htb!] 
\includegraphics[height = 6cm, width = 
\textwidth]{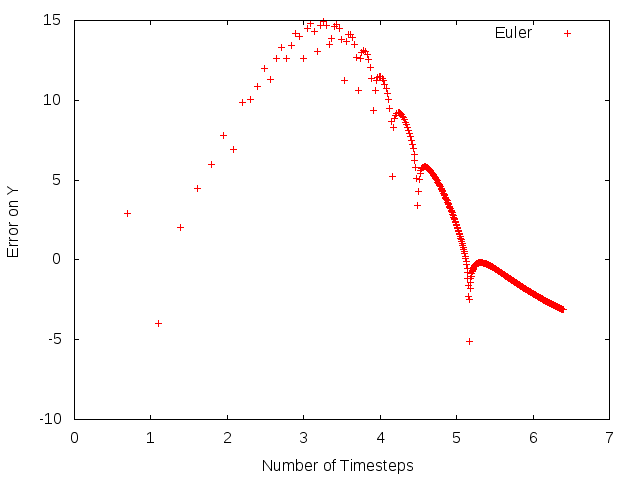} %[width = \textwidth]
\caption{Euler Scheme, $(b=5,T=10)$} \label{fig intro cz5T10N600}
\end{figure}

On Figure \ref{fig intro cz1T10N100} and \ref{fig intro cz5T1N100} appears clearly a first transient state and then the asymptotic steady state,
after a number of time step (between $25$ and $35$). This means that the linear convergence is obtained for $h$ being smaller that some $h^*$. This phenomenon appears here if $b$ is 'big' or $T$ is 'big'. On Figure \ref{fig intro cz5T10N600}, where both $b$ and $T$ are 'big', things are even worse. The correct behaviour of the scheme is observed only for very large $n$ ($n$ is larger than $240$).

%\textcolor{red}{$h^*=?$ ce serait bien d'avoir le meme que dans la figure 2}.

We see that even in the 'pure' BSDE setting, i.e. when $f$ depends only on $z$, the numerical method exhibits some instability. In the sequel, we  investigate this unstable behaviour and provide sufficient and necessary (in some sense) conditions in order to avoid it.

%\newpage
\subsection{Main assumptions and stability of BSDEs} 
Before stating a  precise definition of numerical stability, we recall some known sufficient conditions to obtain a bounded solution $\cY$ to \eqref{eq BSDE}. 
Since we are  interested in the numerical behavior of discretization schemes for BSDEs, we simplify our framework by assuming that the generator is deterministic and does not depend on time $t$ (denoted $(y,z) \mapsto f(y,z)$ by an abuse of notation). We also suppose that $f(0,0)=0$.

In the sequel, we shall make use of the following assumptions.% assumptions.

\HYP{fLy}: The function $f$ is Lipschitz continuous with respect to $y$ with Lipschitz constant $L^Y \geq 0$, i.e.
\begin{align*}
 |f(y',z)-f(y,z)|\le L^Y|y-y'|.
\end{align*}

\HYP{fLz}: The function $f$ is Lipschitz continuous with respect to $z$ with Lipschitz constant $L^Z \geq 0$, i.e.
\begin{align*}
 |f(y,z')-f(y,z)|\le L^Z|z-z'|.
\end{align*}

\HYP{fmy}: The function $f$ is monotone in $y$ with a constant of monotonicity $l^Y \geq 0$, i.e.
\begin{align*}
\langle y-y',f(y,z) - f(y',z) \rangle \leq -l^Y \abs{y-y'}^2.
\end{align*}

Moreover $f$ is continuous in $y$ and has a controlled growth in $y$, precisely there exists an increasing function $\kappa : \R^+ \rightarrow \R^+$ such that
\begin{align*}
 \abs{f(y,z)} \leq \abs{f(0,z)} + \kappa(\abs{y}).
\end{align*}

Let us remark that if assumptions \HYP{fLy} and \HYP{fmy} hold true then we have $l^Y \leq L^Y$.

%\textcolor{blue}
{
\begin{Proposition}[Stability of BSDEs] \label{pr estim BSDE} 
Let us assume that \HYP{fLz}-\HYP{fmy} hold true and $\|\xi \|_\infty<\infty$, where $\|\cdot\|_\infty$ denotes the $L^\infty$ norm of the euclidian norm of random vector.
 Then, if $m=1$, there exists a unique solution $(\cY,\cZ) \in \SP{2} \times \HP{2}$ such that 
 \begin{align}
 \label{estimee solution multidim}
 \|\cY\|_{\SP{\infty}}:= \textrm{essup}_{0 \leq t \leq T} \abs{\cY_t} \leq \|\xi\|_{\infty} \;,
\end{align}
and if $m>1$, the previous statement holds true assuming moreover that 
\begin{equation*}
 (L^Z)^2 \leq 2 l^Y.
\end{equation*}
\end{Proposition}
}

\begin{Remark}
 It is possible to obtain the same type of result as Proposition \ref{pr estim BSDE}  when $f$ is not Lipschitz with respect to $z$ but has quadratic growth (see e.g. \cite{Kobylanski-00,Lepeltier-SanMartin-98}).We do not discuss here this result because when computing numerical approximations of quadratic BSDEs, a first step consists in truncating the generator with respect to $z$ to obtain a Lipschitz generator (see e.g. \cite{Richou-12,Chassagneux-Richou-14}).
\end{Remark}

\paragraph{Proof of Proposition \ref{pr estim BSDE}}
For this proposition the existence and the uniqueness of the solution come from \cite{Pardoux-99}. The estimate \eqref{estimee solution multidim} is quite standard to obtain since it is just sufficient to apply the Ito formula to the process $e^{ ((L^Z)^2 -2 l^Y)t}\abs{\cY_t}^2$ (see e.g. Proposition 2.2 in \cite{Pardoux-99}). In dimension $m=1$, the estimate \eqref{estimee solution multidim} comes from a classical linearization argument: see e.g. \cite{Briand-Hu-98} or \cite{Royer-04}.
\cqfd

\vspace{2mm}

From now on, we assume that \HYP{fLz} and \HYP{fmy} hold true, $\|\xi\|_{\infty}<+\infty$ and 
\begin{equation}
 \label{hypothese coefficients multidim}
 (L^Z)^2 \leq 2 l^Y,
\end{equation}
when $m>1$. Then, from Proposition \ref{pr estim BSDE}, we have  $\|\cY\|_{\SP{\infty}} \leq \|\xi\|_{\infty}$.

\vspace{2mm}

\subsection{Definition: Numerical Stability}

In practice, we will study the numerical stability of the following family of schemes.
\\
For $n\ge 1$, we set $\pi = \set{t_0:=0,\dots,t_i,\dots,t_n := T}$ a discrete-time grid
of $[0,T]$. We denote $h_i := \ti{+1}-\ti{}$ and  $\max h_i = h$. We assume that $h=O(\frac1n)$. On a probability space $(\wh{\Omega},\wh{\cA},\wh{\P})$, we are given discrete-time filtration
$\wh{\F}=(\wh{\cF}_{\ti{}})_{1\le i \le n}$ associated to $\pi$.

\begin{Definition}  \label{de scheme} 

(i) The terminal condition of the scheme is given by $\wh{\xi}$ which is an $\wh{\cF}_T$-measurable square integrable random variable. 

(ii) The transition from step $i+1$ to step $i$ is given by
\begin{align*}
{Y}_i &= \whEFp{\ti{}}{{Y}_{i+1}^{} +h_i \theta f({Y}_i^{},{Z}_i^{}) + h_i (1-\theta)  f({Y}_{i+1}^{},{Z}_i^{}) }
\\
{Z}_i^{} &= \whEFp{\ti{}}{  {Y}_{i+1}^{}{H}_i' }.
\end{align*}
for $\theta \in \set{0,1}$. The value $\theta=1$ corresponds to the implicit scheme and $\theta=0$ to the ``pseudo-explicit'' scheme.

We assume that $H$-coefficients $({H}_i)_{0 \leqslant i <n}$ are some $\R^{d}$  independent random vectors such that, for all $0 \leqslant i < n$, ${H}_i$ is $\wh{\cF}_{t_{i+1}}$ measurable, $ \whEFp{\ti{}}{{H}_i{}}=0$,
 \begin{align}
 c_i I_{d \times d}= h_i \whEFp{}{{H}_i  {H}_i'}= h_i \whEFp{\ti{}}{{H}_i {H}_i'}\,, \label{eq ass H 2}
 \end{align} 
and 
 \begin{align}
 \frac{\lambda}{d} \leqslant c_i \leqslant \frac{\Lambda}{d} \,, \label{eq ass H 3}
 \end{align} 
 where $\lambda$, $\Lambda$ are positive constants which do not depend on $T$ and $n$. Let us remark that \eqref{eq ass H 2} and \eqref{eq ass H 3} imply that
\begin{align}
 \lambda \leqslant h_i \whEFp{}{|{H}_i|^2}= h_i \whEFp{\ti{}}{|{H}_i|^2}\leqslant \Lambda\,. \label{eq ass H}
 \end{align}

\end{Definition}

%\textcolor{blue}
{
We would like to discuss now the well-posedness of the above methods, meaning:
\begin{itemize}
\item one can solve for ${Y}_i$ in practice, when the scheme has an implicit feature i.e. when $\theta=1$;
\item for all $i \le n-1$, $({Y}_i,{Z}_i)$ are square-integrable.
\end{itemize}
Note that, in the sequel, we will always assume the well-posedness of schemes given in Definition \ref{de scheme}.

In the following lemma, we recall sufficient conditions to obtain this property.

{
\begin{Lemma} \label{le scheme wellposedness}
\begin{enumerate}[(i)]
\item For $\theta=0$ (explicit scheme), the scheme is well-posed under \HYP{fLy}-\HYP{fLz}.
\item For $\theta=1$ (implicit scheme), the scheme is well-posed under \HYP{fmy}-\HYP{fLz}.
\end{enumerate}
\end{Lemma}
}

%\textcolor{red}
{
\proof
Statement (i) follows directly. Statement (ii) is more involved. Let us assume that the scheme is well-posed until step $i+1$ and let us show that $({Y}_i, {Z}_i)$ are well-defined and square integrable. Obviously, there is no issue for ${Z}_i$. By remarking that the map $F : y \mapsto y-\theta f(y,{Z}_i(\omega))$ is almost surely strongly monotone since we have
$$\langle y'-y,F(y')-F(y) \rangle \geq (1+h_i l^Y)\abs{y'-y}^2, \quad \forall y,y' \in \R^m,$$
we can use same arguments as in section 4.4 of \cite{Lionnet-dosReis-Szpruch-14} to show the existence of a unique $\widehat{\mathcal{F}}_{t_i}$-measurable r.v. ${Y}_i$ such that 
$${Y}_i = \whEFp{\ti{}}{{Y}_{i+1}^{} +h_i f({Y}_i^{},{Z}_i^{})}.$$
Moreover, we have
\begin{align*}
 \abs{{Y}_i}^2 &= \whEFp{\ti{}}{{Y}_{i}'{Y}_{i+1} +h_i {Y}_{i}'f({Y}_i^{},{Z}_i^{})}\\
 &\leq \frac{\abs{{Y}_i}^2}{2} +  (L^Z)^2h^2\abs{{Z}_i}^2+\whEFp{\ti{}}{\abs{{Y}_{i+1}^{}}^2}
\end{align*}
and so ${Y}_i$ is a square integrable r.v. 
\cqfd
}

\begin{Remark}
\begin{enumerate}[(i)]
\item 
One could extend Definition \ref{de scheme} to any $\theta \in [0,1]$ and actually carry on the analysis made in the next section. One should note however that this would not be the usual $\theta$-scheme as only $Z_i$ appears in the approximation. In particular, one cannot hope to retrieve
an order 2 scheme for $\theta=\frac12$ as in \cite{criman14}, in the general case where $f$ depends on $Z$.
\item The theoretical discrete-time approximation of \eqref{eq th scheme Y imp}-\eqref{eq th scheme Y exp} and \eqref{eq th scheme Z} belong to the above class of approximations. It suffices to work on $(\Omega,\cA,\P)$,
to set $\widehat{\cF}_{\ti{}}=\cF_{\ti{}}$ and $H_i{}=\frac{W_{\ti{+1}}-W_{\ti{}}}{\ti{+1}-\ti{}}$.

\item The above setting encompasses the case of tree methods like cubature methods, see e.g. \cite{criman12}, and quantization methods, see e.g. \cite{Bally-Pages-03}.
\end{enumerate}
\end{Remark}

We now introduce the notion of numerical stability which is an attempt to formalise the  phenomenon described in section \ref{subse motivation}. Roughly speaking, for an ODE we say that a scheme is numericaly stable if the numerical solution obtained with this scheme remains bounded when the real solution is bounded. It is not possible to transpose directly this notion to BSDEs since here $T$ is fixed but we defined a closely related notion of numerical stability. To do this, let us consider a BSDE such that $Y$ is bounded by a constant that does not depend on $T$, recalling Proposition \ref{pr estim BSDE}. Roughly speaking, we will say that a scheme is numericaly stable if the numerical approximation of this BSDE remains bounded by a bound that does not depend on $T$ as well.

\begin{Definition}[Numerical Stability] 
\label{de num stab}
We say that the scheme given in Definition \ref{de scheme} is numerically stable if, there exists $h^*$ such that
for all $h \le h^*$,
\begin{align*}
|{Y}_0|\le \|\widehat{\xi}\|_\infty\;,
\end{align*}
for all essentially bounded $\hat{\cF}_T$-measurable random variable $\widehat{\xi}$.
\end{Definition}

We  also introduce an unconditional stability property for the scheme above.

\begin{Definition}[A-stability] \label{de A stab}
We say that the scheme is $A$-stable, if $h^*=\infty$ in the definition above.
\end{Definition}

\begin{Remark}
 If $\xi$ is non random, the schemes given in Definition \ref{de scheme} are the usual implicit ($\theta=1$) and explicit ($\theta=0$) Euler schemes for ODEs.
 Results for numerical stability are well known, see e.g. \cite{Butcher-03}. In particular, the implicit Euler scheme is A-stable and the explicit Euler scheme is stable if $\frac{|L^Y|^2}{2l^Y} h \le 1$. We will show that the so-called 'implicit' Euler scheme for BSDEs may not be A-stable.
\end{Remark}

\begin{Remark}
Let us mention that a notion of $L^2$-stability has already been introduced for the above method \cite{Chassagneux-Crisan-13} and extended in \cite{Lionnet-dosReis-Szpruch-14}. This notion does not coincide with the one considered here. 
Indeed, it allows only to prove convergence of the scheme, focusing on the asymptotic $h \rightarrow 0$.
\end{Remark}

%\input{lineares}

%\input{BTZstability}

%!TEX root = main.tex 

\section{Sufficient conditions for numerical stability}

We present here our main results concerning the numerical stability of the methods given in Definition \ref{de scheme}.
The conditions below allow to determine the range of timesteps $h>0$ for which the methods are guaranteed to be stable.
We state our results by considering separately the multidimensional setting and the one-dimensional setting for $Y$.
Similarly to the continuous BSDEs case, we obtain stability results using different sets of assumption.
In the next Section, we will perform a Von Neumann stability analysis, which completes the results of this section.

\subsection{Multidimensional case}

In this paragraph and the next one, we assume that the scheme given in Definition \ref{de scheme} is well-posed, see Lemma \ref{le scheme wellposedness} for sufficient conditions. Our first result concerns the multidimensional case for $Y$.

\begin{Proposition}
\label{prop cond suff multidim}
 Assume that \HYP{fLz} and \HYP{fmy} hold with $l^Y>0$ and if $\theta =0$, that \HYP{fLy} is in force as well. 
 %Moreover, we assume that assumption \HYP{fLy} hold when . 
If, moreover, 
 \begin{equation}
  \label{condition suffisante multidim}
  \frac{\left(\sqrt{\Lambda} L^Z + \sqrt{h} L^Y (1-\theta)\right)^2}{2l^Y} \leq 1,
 \end{equation}
then the scheme given in Definition \ref{de scheme} is numericaly stable, recalling Definition \ref{de num stab}.
\end{Proposition}

Before giving the proof of the above proposition, we make the following observations.
\begin{Remark}
 \begin{enumerate}[(i)]
  \item The best sufficient condition is obtained for the implicit scheme ($\theta=1$). In this case, \eqref{condition suffisante multidim} becomes 
  \begin{equation}
  \label{condition suffisante multidim pour theta=1}
  \Lambda (L^Z)^2 \leqslant 2 l^Y
  \end{equation}
  which is exactly the assumption \eqref{hypothese coefficients multidim} when $\Lambda =1$. Moreover, \eqref{condition suffisante multidim pour theta=1} does not depend on $h$ which means that when this condition is satisfied then the scheme is A-stable, recalling Definition \ref{de A stab}.
  \item The fact that we do not need assumption \HYP{fLy} when $\theta=1$ (implicit scheme) allows us to study the stability of the untruncated implicit scheme for BSDEs with polynomial growth drivers with respect to $y$ introduced and studied in \cite{Lionnet-dosReis-Szpruch-14}.  
 \end{enumerate}
\end{Remark}

\paragraph{Proof of Proposition \ref{prop cond suff multidim}}
For $0 \le i \le n-1$, setting 
\begin{equation*}
 \Gamma_i :=  \E_{t_i}\left[\left\{Y_{i+1} + h_i(1-\theta)f(Y_{i+1},Z_i)\right\}H_i'\right],
\end{equation*}
As in the seminal paper \cite{Briand-Delyon-Memin-02}, we observe that
\begin{equation*}
Y_i = Y_{i+1} + h_i \left( \theta f(Y_i,Z_i) + (1-\theta) f(Y_{i+1},Z_i) \right) - h_i c_i^{-1} \Gamma_iH_i - \Delta M_i,
\end{equation*}
with $\E_{t_i}\left[\Delta M_i\right] = 0$ and $\E_{t_i}\left[\Delta M_iH_i'\right] = 0$.
Using the identity $|y|^2 = |x|^2 + 2x'(y-x) + |y-x|^2$ with $y= Y_{i+1}$ and $x = Y_i$, and taking expectation on
both sides, we compute
\begin{eqnarray}
 \nonumber \abs{Y_i}^2 &=& \E_{t_i} \left[\abs{Y_{i+1}}^2 + 2 h_i Y_i'\left\{ \theta f(Y_i,Z_i)+(1-\theta) f(Y_{i+1},Z_i) \right\} -\abs{Y_i -Y_{i+1}}^2 \right]\\
 \nonumber &=& \E_{t_i} \left[\abs{Y_{i+1}}^2 + 2 h_i Y_i' f(Y_i,0) +2 h_i Y_i' \left\{f(Y_i,Z_i)-f(Y_i,0) \right\}\right.\\
 \label{equation ito discret carre} && \quad \quad \left.+2(1-\theta)h_iY_i'\left\{f(Y_{i+1},Z_i)-f(Y_i,Z_i)\right\} -\abs{Y_i -Y_{i+1}}^2 \right].
\end{eqnarray}
Then assumptions \HYP{fLz}, \HYP{fmy} and \HYP{fLy} (if $\theta=0$) on $f$ yield 
\begin{eqnarray*}
 \abs{Y_i}^2 &\leq& \E_{t_i} \left[\abs{Y_{i+1}}^2 -2l^{Y} h_i \abs{Y_i}^2 + 2 h_i L^Z \abs{Y_i}\abs{Z_i}+2h_i (1-\theta) L^Y \abs{Y_i}\abs{Y_{i+1}-Y_i}\right.\\
 && \quad \quad \left. -\abs{Y_{i+1}-Y_i}^2 \right].
\end{eqnarray*}
We now introduce two constants $\alpha >0$ and $\beta>0$ to be set latter on. Using Young inequality (twice), we compute
\begin{eqnarray*}
 \abs{Y_i}^2 &\leq& \E_{t_i} \left[\abs{Y_{i+1}}^2 +(\alpha+\beta-2l^{Y}) h_i \abs{Y_i}^2 + \frac{h_i (L^Z)^2}{\alpha}  \abs{Z_i}^2+ \frac{h_i (1-\theta)^2 (L^Y)^2}{\beta} \abs{Y_{i+1}-Y_i}^2\right.\\
 && \quad \quad \left. -\abs{Y_{i+1}-Y_i}^2 \right].
\end{eqnarray*}
Let us remark that if $\theta=1$ or $L^Y=0$ we do not need to introduce $\beta$ and to use the second Young inequality. In the same way, if $L^Z=0$, we do not need to introduce $\alpha$ and to use the first Young inequality. Since  $Z_i = \E_{t_i} \left[(Y_{i+1}-Y_i) H_i' \right]$, we apply Cauchy-Schwarz inequality to obtain
\begin{equation}
 h_i \abs{Z_i}^2 \leqslant \Lambda \E_{t_i} \left[ \abs{Y_{i+1}-Y_i}^2\right]
\end{equation}
and then the previous inequality becomes
\begin{eqnarray*}
 \abs{Y_i}^2 &\leq& \E_{t_i} \left[\abs{Y_{i+1}}^2 +(\alpha+\beta-2l^{Y}) h_i \abs{Y_i}^2 + \left(\frac{\Lambda (L^Z)^2}{\alpha}+ \frac{h (1-\theta)^2 (L^Y)^2}{\beta}-1\right) \abs{Y_{i+1}-Y_i}^2\right].
\end{eqnarray*}
Finally setting in the above inequality
\begin{eqnarray*}
 \alpha &=& \Lambda (L^Z)^2+ \sqrt{\Lambda h(L^Z)^2(1- \theta)^2(L^Y)^2},\\
 \beta &=& h (1-\theta)^2 (L^Y)^2+ \sqrt{\Lambda h(L^Z)^2(1- \theta)^2(L^Y)^2},
\end{eqnarray*}
leads to
\begin{eqnarray*}
 \abs{Y_i}^2 &\leq& \E_{t_i} \left[\abs{Y_{i+1}}^2 +\left(\left(\sqrt{\Lambda}L^Z+\sqrt{h}L^Y(1-\theta)\right)^2-2l^{Y}\right) h_i \abs{Y_i}^2\right].
\end{eqnarray*}
Under assumption \eqref{condition suffisante multidim} we get 
\begin{equation*}
 \abs{Y_i}^2 \leq \E_{t_i} \left[\abs{Y_{i+1}}^2 \right]
\end{equation*}
and an easy induction concludes the proof.
\cqfd

\subsection{One-dimensional case}
We now turn to the one-dimensional setting for $Y$ and prove the numerical stability of the Euler scheme under slightly different assumptions.

\begin{Proposition}
\label{prop cond suff unidim}
 Assume that $m=1$ and assumptions \HYP{fLz} and \HYP{fmy} hold true. Moreover, when $\theta =0$ we assume that assumption \HYP{fLy} holds and $l^Y>0$. Finally, we also suppose that
  \begin{equation}
  \label{condition suffisante unidim}
  h\left[ \frac{(1-\theta)^2(L^Y)^2}{2l^Y} + L^Z (\max_{0 \leq i \leq n-1} |H_i| )\right] \leq 1.
 \end{equation}
 Then the scheme given in Definition \ref{de scheme} is numerically stable, recalling Definition \ref{de num stab}.
\end{Proposition}

\begin{Remark}
 \begin{enumerate}[(i)]
  \item Assumption \eqref{condition suffisante unidim} imposes that $H$ is bounded.
  \item The best sufficient condition is obtained for the implicit scheme ($\theta=1$). Nevertheless, even in this case,  condition \eqref{condition suffisante unidim}
  does not guarantee $A$-stability, recall Definition \eqref{de A stab}.

  \item When $\theta=1$ (implicit scheme), a comparison theorem holds: see Proposition 2.4 and Corollary 2.5 in \cite{Chassagneux-Richou-14}.
  \item It is worth to compare condition \eqref{condition suffisante unidim} to \eqref{condition suffisante multidim}. First of all, when $f$ does not depend on $z$, $L^Z=0$ and then assumptions \eqref{condition suffisante unidim} and \eqref{condition suffisante multidim} are equal: we find the classical stability condition for ODEs, that is to say 
  $$h\frac{(1-\theta)^2 (L^Y)^2}{2l^Y} \leq 1.$$

 In the general case, it is important to remark that condition \eqref{condition suffisante unidim} is fulfilled as soon as $h$ is small enough whereas it is not the case for \eqref{condition suffisante multidim}. 
  
 \item Since assumption \HYP{fLy} is not required when $\theta=1$ (implicit scheme), our result can be applied to study the stability of the untruncated implicit scheme for BSDEs with polynomial growth drivers with respect to $y$ introduced in \cite{Lionnet-dosReis-Szpruch-14}. 

 \end{enumerate}
\end{Remark}

\paragraph{Proof of Proposition \ref{prop cond suff unidim}}
We  adapt the proof of Proposition \ref{prop cond suff multidim} to the one-dimensional setting. Let us denote, for $0 \leq i \leq n-1$, 
\begin{equation*}
 \gamma_i = \frac{f(Y_i,Z_i)-f(Y_i,0)}{\abs{Z_i}^2}Z_i \mathbbm{1}_{\set{Z_i \neq 0}}.
\end{equation*}
Then, using the definition of $Z_i$, equality \eqref{equation ito discret carre} becomes
\begin{eqnarray*}
 \abs{Y_i}^2  &=& \E_{t_i} \left[\abs{Y_{i+1}}^2 + 2 h_i Y_i f(Y_i,0) +2 h_i Y_i Y_{i+1}\gamma_i H_i \right.\\
  && \quad \quad \left.+2(1-\theta)h_iY_i\left\{f(Y_{i+1},Z_i)-f(Y_i,Z_i)\right\} -\abs{Y_i -Y_{i+1}}^2 \right].
\end{eqnarray*}
Observing that 
$$2Y_i(Y_{i+1}-Y_i)=\abs{Y_{i+1}}^2-\abs{Y_{i+1}-Y_i}^2-\abs{Y_i}^2$$
and $\E_{t_{i}}[\abs{Y_i}^2 \gamma_i H_i]=0$, we compute
\begin{eqnarray*}
 \abs{Y_i}^2  &=& \E_{t_i} \left[(1+h_i\gamma_i H_i)\abs{Y_{i+1}}^2 + 2 h_i Y_i f(Y_i,0) +2(1-\theta)h_iY_i\left\{f(Y_{i+1},Z_i)-f(Y_i,Z_i)\right\}\right.\\
  && \quad \quad \left.-(1+h_i\gamma_iH_i)\abs{Y_i -Y_{i+1}}^2 \right].
\end{eqnarray*}
Using assumptions \HYP{fLz}, \HYP{fmy} and \HYP{fLy} on $f$ together with Young inequality, we get, for all $\alpha>0$,
\begin{eqnarray*}
 \abs{Y_i}^2  &\leq& \E_{t_i} \left[(1+h_i\gamma_i H_i)\abs{Y_{i+1}}^2 + h_i\left( \frac{(1-\theta)^2(L^Y)^2}{\alpha}-2l^Y \right)\abs{Y_i}^2\right.\\
  && \quad \quad \left.-(1+h_i\gamma_iH_i-\alpha h_i)\abs{Y_i -Y_{i+1}}^2 \right].
\end{eqnarray*}
Let us remark that we do not need to introduce $\alpha$ and  use Young inequality if $\theta=1$ or $L^Y=0$. Otherwise, setting $\alpha = (1-\theta)^2(L^Y)^2/(2l^Y)$, we obtain
\begin{eqnarray*}
 \abs{Y_i}^2  &\leq& \E_{t_i} \left[(1+h_i\gamma_i H_i)\abs{Y_{i+1}}^2 -\left(1+h_i\gamma_iH_i- h_i\frac{(1-\theta)^2(L^Y)^2}{2l^Y}\right)\abs{Y_i -Y_{i+1}}^2 \right].
\end{eqnarray*}
Since we assume that \eqref{condition suffisante unidim} and assumption \HYP{fLz} on $f$ hold, then we have 
\begin{equation*}
 1+h_i\gamma_iH_i- h_i\frac{(1-\theta)^2(L^Y)^2}{2l^Y} \geq 1-h_iL^Z\abs{H_i}- h_i\frac{(1-\theta)^2(L^Y)^2}{2l^Y} \geq 0,
\end{equation*}
and
\begin{eqnarray*}
 \abs{Y_i}^2  &\leq& \E_{t_i} \left[(1+h_i\gamma_i H_i)\abs{Y_{i+1}}^2 \right].
\end{eqnarray*}
An easy induction leads to 
\begin{align*}
|Y_i|^2 &\leq \E_{t_i} \left[ \prod_{j=i}^{n-1} (1+ h_j \gamma_j H_j) |Y_n|^2\right].
%\\
%& \le \EFp{\ti{}}{ \Prod_{j=i}^{n-1} (1+ h_j \gamma_j H_j)} ||\,|Y_n|^2||_\infty
\end{align*}
Finally, for all $0 \leq i \leq n-1$,  \eqref{condition suffisante unidim} and assumption \HYP{fLz} on $f$ yields that 
\begin{equation*}
 1+h_i\gamma_i H_i \geq 1- h_iL^Z\abs{H_i} \geq 0.
\end{equation*}
We then easily obtain the inequality
\begin{align*}
|Y_i|^2 &\leq \E_{t_i} \left[ \prod_{j=i}^{n-1} (1+ h_j \gamma_j H_j) \right] \|Y_n\|_{\infty}^2 \leq  \|Y_n\|_{\infty}^2,
\end{align*}
proving the numerical stability of the scheme.
\cqfd

%\input{vonneum}

%!TEX root = main.tex 

\section{Von Neumann stability analysis}
\label{se VN}

\newcommand{\ic}{\mathbf{i}}

In this section, we will perform a Von Neumann stability analysis, inspired by what is done for PDE.
We will restrict our study to the one dimensional case for $Y$ i.e. $m=1$. Moreover, we shall assume here  a uniform time step: $h_i=h$ for all $0 \leq i <n$. The analysis is performed by considering $\mathbb{C}$-valued terminal conditions as explained below.
We thus work in the setting of the previous sections extended to one-dimensional $\mathbb{C}$-valued BSDEs.

We now define the Von Neumann stability for BSDE schemes in our framework. 

\begin{Definition}[Von Neumann Stability] \label{de vn stab}
For $h>0$ we say that the scheme given in Definition \ref{de scheme} is Von Neumann stable (also denoted VN stable) if for all $k \in \mathbb{R}^d$, we have $\abs{Y_0} \leq 1$ when $\xi:=e^{\ic \sum_{\ell=1}^d k_{\ell} W^{\ell}_T}$. We call VN stability region the set of all $h>0$ such that the scheme is VN stable. If this VN stability region is equal to $]0,+\infty[$ then we say that the scheme is VN A-stable.
\end{Definition}
It is clear that numerical stability previously studied implies VN stability, once extended to $\mathbb{C}$-valued BSDEs. In this section, we will perform the Von Neumann stability analysis considering only linear mapping $f$ i.e.
\begin{align*}
 f(y,z) = ay + \sum_{\ell=1}^d b_\ell z^\ell,
\end{align*}
with $a\leq 0$ and $b \in \mathbb{R}^d$. 

Moreover we will only study the classical  scheme given in Definition \ref{de 
scheme} with $H_i=h^{-1} (W_{t_{i+1}}-W_{t_i})$, 
$(\wh{\Omega},\wh{\cA},\wh{\P})=(\Omega,\cA,\P)$ and $\wh{\cF}=\cF$. 
%and $\theta \in \set{0,1}$, i.e. the implicit and pseudo-explicit schemes. 

We observe then that the $(H_i)_{0 \leq i <n}$ are unbounded, so the unidimensional sufficient condition \eqref{condition suffisante unidim} cannot be fulfilled. On the other hand, the multidimensional sufficient condition \eqref{condition suffisante multidim} becomes, in our framework,
\begin{equation}
 \label{condition suffisante multidim cas lineaire}
 \frac{\left(\textcolor{black}{\sqrt{d}}\abs{b} + \sqrt{h} \abs{a} (1-\theta)\right)^2}{2\abs{a}} \leq 1.
\end{equation}
Obviously, \eqref{condition suffisante multidim cas lineaire} is a too strong assumption in practice for the unidimensional case. The VN stability analysis performed below allows us to identify necessary conditions for the numerical stability of the Euler scheme. Importantly, we shall observe that those conditions depend on the dimension of the $Z$ process, even in the one-dimensional case for $Y$.
%, which has important consequences when someone uses these numerical schemes for solving high dimensional problems.

\subsection{Von Neumann stability analysis of the implicit Euler scheme}
We study here the implicit Euler scheme i.e. the scheme given in Definition \ref{de scheme} with $\theta=1$.
 Let us define
\begin{equation}
 \label{definition norme infinie b}
 \abs{b}_{\infty} = \max \left( \abs{b^+}, \abs{b^-}\right),
\end{equation}
with $b^+=(b_1 \vee 0,...,b_d \vee 0)$ and $b^-=(b_1 \wedge 0,...,b_d \wedge 0)$. 

We then have the following results concerning the VN stability of the implicit Euler scheme.
\begin{Proposition}
\label{prop critere stabilite VN implicite}
 \begin{enumerate}[(i)]
  \item If $\abs{b}=0$, then the scheme is VN A-stable.
  \item Assume that $\abs{b}>0$. Then the scheme is VN stable if and only if $\abs{b}_{\infty}^2 h \leq 1$ or,  $\abs{b}_{\infty}^2 h > 1$ and
  \begin{equation}
  \label{critere stabilite VN implicite}
   (1-ah)^2-\abs{b}_{\infty}^2 h e^{\frac{1}{\abs{b}_{\infty}^2 h}-1} \geq 0.
  \end{equation}
  \item In particular, when $a=0$, i.e. when $f$ only depends on $z$, the scheme is VN Stable if and only if $\abs{b}_{\infty}^2 h \leq 1$.
 \end{enumerate}
\end{Proposition}

\paragraph{Proof of Proposition \ref{prop critere stabilite VN implicite}}
Using an induction argument, we first show that  $Y_i = y_i e^{\ic \sum_{\ell=1}^d k_{\ell} W^{\ell}_{t_i}}$ with $y_i \in \mathbb{C}$. Indeed, we observe that this is true for $Y_n$ with $y_n=1$, recalling Definition \ref{de vn stab}. Then, if  $Y_{i+1} = y_{i+1} e^{\ic \sum_{\ell=1}^d k_{\ell} W^{\ell}_{t_{i+1}}}$, we compute
\begin{align*}
(1-ah)Y_i &= \mathbb{E}_{t_i} \left[ y_{i+1} \left(1+ \sum_{\ell =1}^d b_{\ell} \Delta W_i^{\ell}\right)e^{\ic \sum_{\ell=1}^d k_{\ell} W^{\ell}_{t_{i+1}}} \right]\\
  &= y_{i+1}\mathbb{E} \left[ \left(1+ \sum_{\ell =1}^d b_{\ell} \Delta W_i^{\ell}\right)e^{\ic \sum_{\ell=1}^d k_{\ell} \Delta W^{\ell}_{i}} \right]e^{\ic \sum_{\ell=1}^d k_{\ell} W^{\ell}_{t_{i}}}\\
  &= y_{i+1} \left(1+ \ic h \sum_{\ell =1}^d b_{\ell} k_\ell \right)e^{-\frac{\sum_{\ell=1}^d k_{\ell}^2 h}{2}} e^{\ic \sum_{\ell=1}^d k_{\ell} W^{\ell}_{t_{i}}}.
\end{align*}
%\textcolor{red}{mettre definition increment mouvement brownien (dans l'intro ?)}
Thus we have $Y_i = y_i e^{\ic \sum_{\ell=1}^d k_{\ell} W^{\ell}_{t_i}}$ with
$$y_i=\frac{\left(1+ \ic h \sum_{\ell =1}^d b_{\ell} k_\ell \right)e^{-\frac{\sum_{\ell=1}^d k_{\ell}^2 h}{2}}}{1-ah}y_{i+1}:=\lambda y_{i+1},$$
and so we obtain $Y_0=\lambda^n$. Recalling Definition \ref{de vn stab}, we get that the scheme is VN stable if and only if, for all $k \in \R^d$ we have $\abs{\lambda}^2 \leq 1$, i.e.
%\begin{equation*}
%e^{-\sum_{\ell =1}^d x_\ell} \left(1+\left(\sum_{\ell =1}^d b_{\ell} \sqrt{x_\ell}\right)^2h\right) \leq (1-ah)^2, \quad \forall (x_1,...,x_d) \in (\R^+)^d,
%\end{equation*}
%or, equivalently,
\begin{equation}
\label{condition generale stab VN schema implicite}
\varphi(x_1,...,x_d):=(1-ah)^2- e^{-\sum_{\ell =1}^d x_\ell} \left(1+\left(\sum_{\ell =1}^d b_{\ell} \sqrt{x_\ell}\right)^2h\right) \geq 0, \quad \forall x \in (\R^+)^d.
\end{equation}
We first  remark that \eqref{condition generale stab VN schema implicite} is always true if $\abs{b}=0$, proving (i). 
We now deal with the case $\abs{b}>0$. Since $\varphi(x) \leq (1-ah)^2$ and $\lim_{\abs{x} \rightarrow +\infty} \varphi(x)=(1-ah)^2$, we know that $\varphi$ is bounded from above and there exists $\tilde{x}$ such that $\inf_{x \in (\R^+)^d} \varphi(x)=\varphi(\tilde{x})$. Then \eqref{condition generale stab VN schema implicite} holds true if and only if $\varphi(\tilde{x}) \geq 0$. We now need to identify $\tilde{x}$.
\begin{enumerate}
 \item If $\tilde{x} \in (\R^{+*})^d$, then $\nabla \varphi(\tilde{x})=0$ and necessarily we must have $b_i/\sqrt{\tilde{x}_i}=b_j/\sqrt{\tilde{x}_j}$ for all $i,j \in \set{1,...,d}$. In particular we must have all $b_i$ with the same sign. If this is true, then $\tilde{x}_i=b_i^2(\frac{1}{\abs{b}^2}-\frac{1}{\abs{b}^4 h})$. Since $\tilde{x} \in (\R^{+*})^d$, this implies that we must have also $\abs{b}^2h >1$. To sum up, if all $b_i$ have the same sign and if $\abs{b}^2h >1$, then $x$ given by ${x}_i=b_i^2\left(\frac{1}{\abs{b}^2}-\frac{1}{\abs{b}^4 h}\right)$ for $1 \leq i \leq d$ is the only candidate for $\tilde{x}$ and then $\varphi(\tilde x)=(1-ah)^2-\abs{b}^2 h e^{\frac{1}{\abs{b}^2 h}-1}$.
 \item We now consider the case where $\tilde x$ is on the boundary of $(\R^{+})^d$. We denote $I$ a non-empty subset of $\set{1,...,d}$ and we assume that $\tilde{x}_i=0$ for $i \in I$ and $\tilde{x}_i>0$ if $i \notin I$. We can use the same reasoning as in the previous step, the only difference coming from the dimension of the space which is strictly smaller. Finally, we obtain that if all $(b_i)_{i \notin I}$ have the same sign and if $h\sum_{i \notin I} b_i^2>1$, then $x$ given by ${x}_i=b_i^2\left(\frac{1}{\abs{b}^2}-\frac{1}{\abs{b}^4 h}\right)\mathbbm{1}_{i \notin I}$ for $1 \leq i \leq d$ is a candidate for $\tilde{x}$ and we have 
 $$\varphi(\tilde x)=(1-ah)^2- h e^{\frac{1}{ h\sum_{i \notin I} b_i^2}-1} \sum_{i \notin I} b_i^2.$$
 \item Since we have a finite number of candidates for $\tilde{x}$, to conclude we just have to compare the value of $\varphi$ for each candidate. Firstly, when $\abs{b}_{\infty}^2 h \leq 1$, then the only candidate is $0$ and so $\varphi(x)\geq \varphi(0)= (1-ah)^2-1 \geq 0$ which implies that the scheme is VN stable. Now let us assume that $\abs{b}_{\infty}^2 h > 1$. By remarking that the function
 $\beta \mapsto (1-ah)^2-\beta h e^{-1+\frac{1}{\beta h}}$ is decreasing on $[1/h,+\infty[$, we obtain that 
 $$\varphi(x) \geq \varphi(\tilde{x})=(1-ah)^2-\abs{b}_{\infty}^2 h e^{\frac{1}{\abs{b}_{\infty}^2 h}-1}.$$ 
\end{enumerate}
This proves $(ii)$ in the statement of the proposition. The remark $(iii)$ can be directly deduced from $(i)$ and $(ii)$ setting $a=0$.
%We treat the case $a=0$ easily thanks to the previous VN stability characterisation.
\cqfd

\vspace{2mm}
We would like now to describe the stability region for $h$. This is easily done for the special case $a=0$ or $b=0$, see $(i)$ and $(iii)$ of the above proposition. 
The following result is a description of the VN stability region in the general case. 
%But it is not straightforward to obtain an explicit VN stability region in the general case. 

\begin{Corollary} \label{co imp VN-stab region}
 There exist real numbers $\tilde{p}>0$ and $\tilde{u}>1$ such that:
 \begin{itemize}
  \item if $-\frac{a}{\abs{b}_{\infty}^2} \geq \tilde{p}$, then \eqref{critere stabilite VN implicite} is true for all $h>0$: the scheme is VN A-stable, recalling Definition \ref{de vn stab};
  \item if $-\frac{a}{\abs{b}_{\infty}^2} < \tilde{p}$, there exists $1<\underline{u} < \tilde{u} < \bar{u}<+\infty$ such that the scheme is VN stable if and only if $h \notin \left]\frac{\underline{u}}{\abs{b}_{\infty}^2},\frac{\bar{u}}{\abs{b}_{\infty}^2}\right[$. Moreover, $\underline{u}$ and $\bar{u}$ are respectively an increasing function and a decreasing function of $p=-\frac{a}{\abs{b}_{\infty}^2}$ satisfying
  $$
  \textcolor{black}{\lim_{-\frac{a}{\abs{b}_{\infty}^2} \rightarrow 0} (\underline{u},\bar{u})=(1,+\infty)}, \quad \lim_{-\frac{a}{\abs{b}_{\infty}^2} \rightarrow \tilde{p}} (\underline{u},\bar{u})=(\tilde{u},\tilde{u}).$$
% \textcolor{red}{$1/b^2??$ ha... il y a deux manieres (au moins...) pour p vers 0... a=0 ou b va a l'infini...}
 \end{itemize}
Numerically we obtain $\tilde{p} \simeq 0.103417$ and $\tilde{u}=7.35491$.\end{Corollary}

\proof
Setting $p=-a/\abs{b}_{\infty}^2$ and $u=\abs{b}_{\infty}^2 h$, then  \eqref{critere stabilite VN implicite}
becomes
\begin{align*}
 \psi(p,u) := (1+pu) - ue^{\frac{1}{u}-1} \ge 0.
\end{align*}
The results are then obtained by studying the sign of the function $u \mapsto \psi(p,u)$ for $u \in (0,\infty)$, 
when $p$ varies in $(0,+\infty)$.  
\cqfd

\vspace{1mm}
The results of Corollary \ref{co imp VN-stab region} are illustrated in Figure \ref{fig imp VN-stab region},  which shows also the point 
$A=(b.\tilde p, \tilde u/b^2)$.

\begin{figure}[htb!] 
\includegraphics[height = 6cm, width = \textwidth]{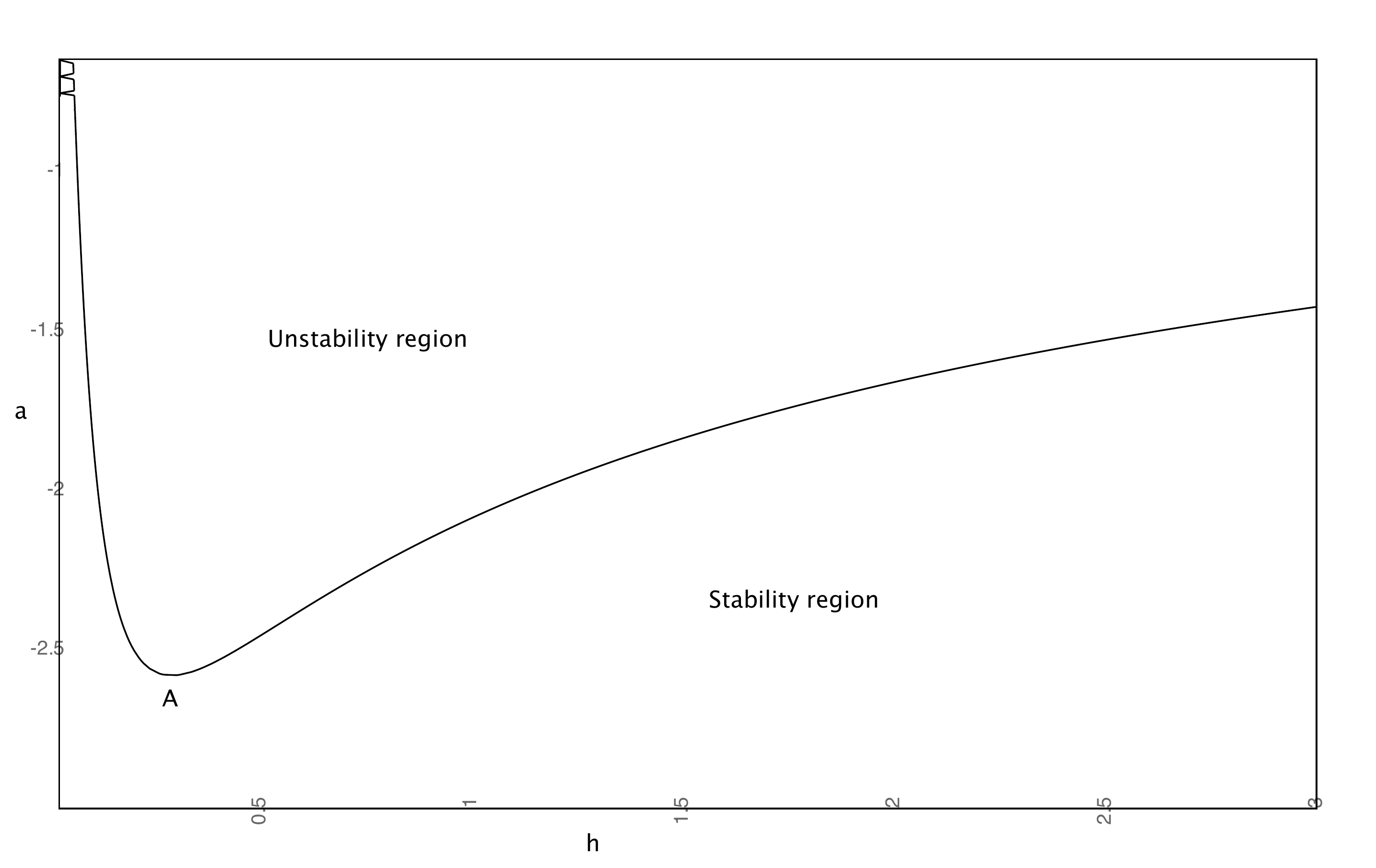} %[width = \textwidth]
\caption{Von Neuman stability region, Implicit Euler, $d=1$, $b=5$} \label{fig imp VN-stab region}
\end{figure}

%\textcolor{magenta}{Add the point $(b.\tilde p, \tilde u/b^2)$ on the figure}

\begin{Remark}
 \begin{enumerate}[(i)]
  \item VN stability regions obtained depend on $\abs{b}_{\infty}$ and so on the dimension $d$. In particular, the VN A-stability condition is more difficult to fulfill when $d$ is large, and the size of the VN stability region decreases with $d$.
  \item The necessary and sufficient VN A-stability condition obtained, namely $-\frac{a}{\abs{b}_{\infty}^2} \geq \tilde{p}$, is much better than the sufficient VN A-stability condition \eqref{condition suffisante multidim cas lineaire}, namely $-\frac{a}{\abs{b}^2} \geq 0.5\textcolor{black}{d}$.
 \end{enumerate}

\end{Remark}

\subsection{Von Neumann stability analysis of the pseudo explicit Euler scheme}

We study here the pseudo explicit Euler scheme i.e. the scheme given in Definition \ref{de scheme} with $\theta=0$.

\begin{Proposition}
\label{prop critere stabilite VN pseudo explicite}
 \begin{enumerate}[(i)]
  \item If $\abs{b}=0$, then the scheme is VN stable if and only if $h \leq -\frac2a$.
  \item Assume that $\abs{b}>0$. Then the scheme is VN stable if and only if $\abs{b}_{\infty}^2 h \leq (1+ah)^2$ and $h \leq -2/a$, or,  $\abs{b}_{\infty}^2 h > (1+ah)^2$ and
  \begin{equation}
  \label{critere stabilite VN pseudo explicite}
   1-\abs{b}_{\infty}^2 h e^{\frac{(1+ah)^2}{\abs{b}_{\infty}^2 h}-1} \geq 0.
  \end{equation}
 \end{enumerate}
\end{Proposition}

\begin{Remark}
  When $a=0$, the same necessary and sufficient condition as in Proposition \ref{prop critere stabilite VN implicite}(iii) holds true, namely the scheme is VN-stable if and only if $\abs{b}_{\infty}^2 h \leq 1$. Indeed, in this case, the implicit scheme and pseudo-explicit scheme are the same.
 \end{Remark}

\paragraph{Proof of Proposition \ref{prop critere stabilite VN pseudo explicite}}
Using the same arguments as in the proof of Proposition \ref{prop critere stabilite VN implicite},  we can write $Y_i = y_i e^{\ic \sum_{\ell=1}^d k_{\ell} W^{\ell}_{t_i}}$ with $y_i \in \mathbb{C}$. Moreover, we compute
\begin{align*}
Y_i &= \mathbb{E}_{t_i} \left[ y_{i+1} \left(1+ ah +\sum_{\ell =1}^d b_{\ell} \Delta W_i^{\ell}\right)e^{\ic \sum_{\ell=1}^d k_{\ell} W^{\ell}_{t_{i+1}}} \right]\\
  &= y_{i+1}\mathbb{E} \left[ \left(1+ah+ \sum_{\ell =1}^d b_{\ell} \Delta W_i^{\ell}\right)e^{\ic \sum_{\ell=1}^d k_{\ell} \Delta W^{\ell}_{i}} \right]e^{\ic \sum_{\ell=1}^d k_{\ell} W^{\ell}_{t_{i}}}\\
  &= y_{i+1} \left(1+ah+ \ic h \sum_{\ell =1}^d b_{\ell} k_\ell \right)e^{-\frac{\sum_{\ell=1}^d k_{\ell}^2 h}{2}} e^{\ic \sum_{\ell=1}^d k_{\ell} W^{\ell}_{t_{i}}}.
\end{align*}
Thus we have $Y_i = y_i e^{\ic \sum_{\ell=1}^d k_{\ell} W^{\ell}_{t_i}}$ with
$$y_i=\left(1+ ah+\ic h \sum_{\ell =1}^d b_{\ell} k_\ell \right)e^{-\frac{\sum_{\ell=1}^d k_{\ell}^2 h}{2}}y_{i+1}:=\lambda y_{i+1},$$
and so we obtain $Y_0=\lambda^n$. Recalling Definition \ref{de vn stab}, we get that the scheme is VN stable if and only if, for all $k \in \R^d$ we have $\abs{\lambda}^2 \leq 1$, i.e.
%\begin{equation*}
%e^{-\sum_{\ell =1}^d x_\ell} \left((1+ah)^2+\left(\sum_{\ell =1}^d b_{\ell} \sqrt{x_\ell}\right)^2h\right) \leq 1, \quad \forall (x_1,...,x_d) \in (\R^+)^d,
%\end{equation*}
%or, equivalently,
\begin{equation}
\label{condition generale stab VN schema pseudo explicite}
\varphi(x_1,...,x_d):=1- e^{-\sum_{\ell =1}^d x_\ell} \left((1+ah)^2+\left(\sum_{\ell =1}^d b_{\ell} \sqrt{x_\ell}\right)^2h\right) \geq 0, \quad \forall x \in (\R^+)^d.
\end{equation}
If $|b|=0$, we observe that\eqref{condition generale stab VN schema pseudo explicite} holds true if and only if $\abs{1+ah} \leq 1$. This proves $(i)$. 
We now study the case $\abs{b}>0$. 
Since $\varphi \leq 1$ and $\lim_{\abs{x} \rightarrow +\infty} \varphi(x)=1$, we have that $\varphi$ is bounded from above and there exists $\tilde{x}$ such that $\inf_{x \in (\R^+)^d} \varphi(x)=\varphi(\tilde{x})$. Then \eqref{condition generale stab VN schema implicite} is true if and only if $\varphi(\tilde{x}) \geq 0$ and it just remains to find $\tilde{x}$. Using the same reasoning as in the implicit case, we finally show that
\begin{enumerate}
 \item if $\abs{b}_{\infty}^2 h \leq (1+ah)^2$ then $\tilde{x}=0$ and so $\varphi(x)\geq \varphi(0)= 1-(1+ah)^2$ which is positive if and only if $h \leq-2/a$,
 \item if $\abs{b}_{\infty}^2 h > (1+ah)^2$ then 
  $$\varphi(x) \geq \varphi(\tilde{x})=1-\abs{b}_{\infty}^2 h e^{\frac{(1+ah)^2}{\abs{b}_{\infty}^2 h}-1}.$$ 
\end{enumerate}
\cqfd

The following Corollary describes the VN stability region more explicitly.

%\textcolor{red}
{
\begin{Corollary}
 \begin{itemize}
  \item If $-\frac{a}{\abs{b}_{\infty}^2} \geq 2$, then the scheme is VN stable if and only if $h \in [0,-2/a]$.
  \item If $-\frac{a}{\abs{b}_{\infty}^2} < 2$, there exists $\overline{h} \in \left[\frac{1}{\abs{b}_{\infty}^2},\frac{-2}{a}\right[ $ such that the scheme is VN stable if and only if $h \in [0,\overline{h}]$. $\overline{h}$ is given by the unique solution of the equation 
  $$1-\abs{b}_{\infty}^2 \overline{h} e^{\frac{(1+a\overline{h})^2}{\abs{b}_{\infty}^2 \overline{h}}-1} = 0.$$
  Moreover, we have
  $$\lim_{-\frac{a}{\abs{b}_{\infty}^2} \rightarrow 0} \overline{h}\abs{b}_{\infty}^2=1, \quad \lim_{-\frac{a}{\abs{b}_{\infty}^2} \rightarrow 2} \overline{h}\frac{\abs{a}}{2}=1.$$
 \end{itemize}
\end{Corollary}
}

\proof
Setting $p=-a/\abs{b}_{\infty}^2$ and $u=\abs{b}_{\infty}^2 h$, the stability region is obtained studying the sign of the function 
$$u \mapsto  1 - ue^{\frac{(1-pu)^2}{u}-1}, \quad u \in \left]\frac{1+2p-\sqrt{1+4p}}{2p^2},\frac{1+2p+\sqrt{1+4p}}{2p^2}\right[,$$
when $p$ varies in $]0,+\infty[$. 
%Indeed, by taking , the previous function becomes \eqref{critere stabilite VN pseudo explicite}. In this way, we obtain the following proposition.
\cqfd

\begin{Remark}
 \begin{enumerate}[(i)]
  \item Once again VN stability regions obtained depend on $\abs{b}_{\infty}$ and so on the dimension $d$. 
  \item Unlike the implicit scheme, the pseudo-explicit scheme is never VN A-stable.
 \end{enumerate}

\end{Remark}

\begin{figure}[htb!] 
\includegraphics[height = 6cm, width = \textwidth]{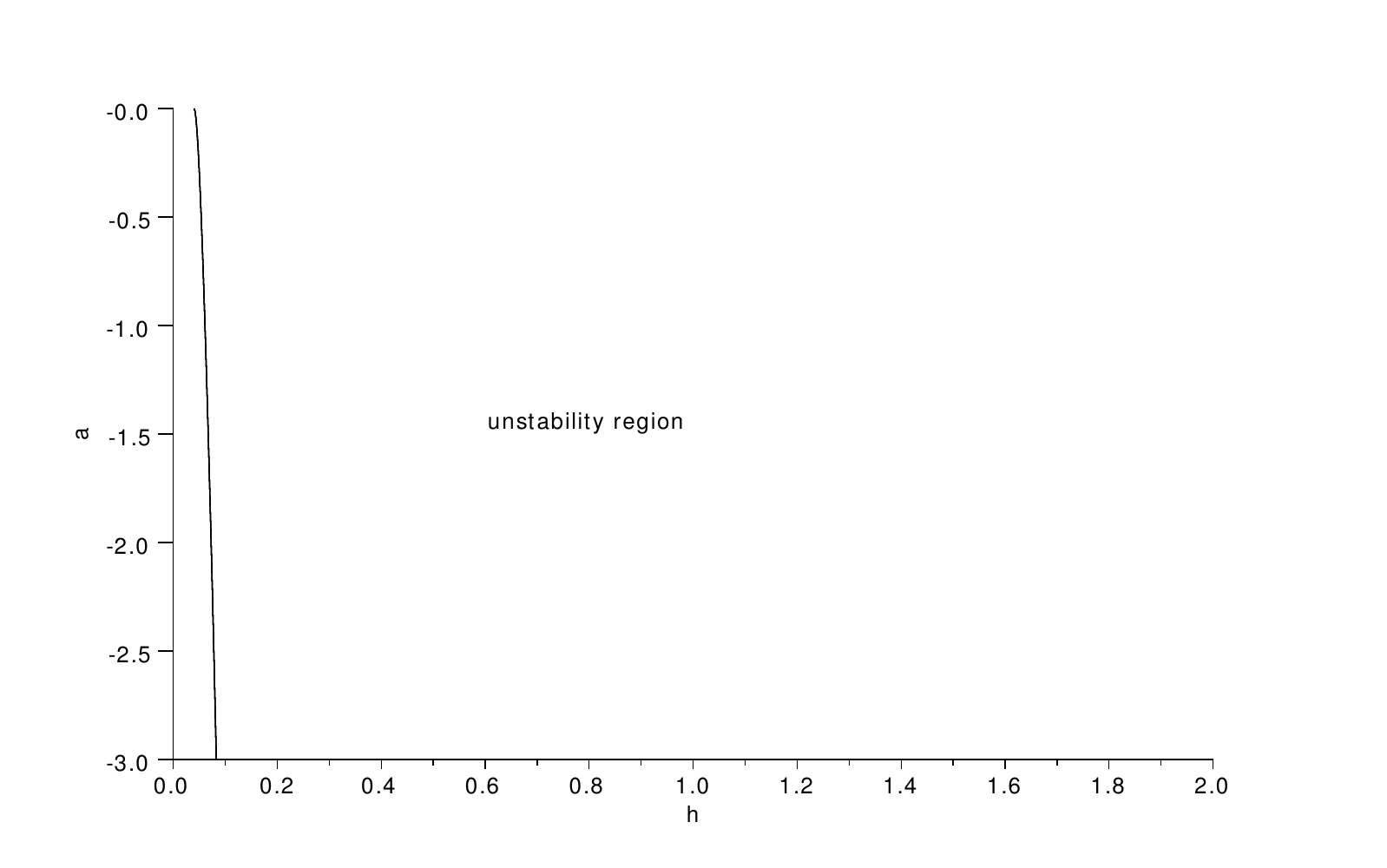} %[width = \textwidth]
\caption{VN stability region, Pseudo-Explicit Euler, $d=1$, $b=5$} \label{fig exp VN-stab region}
\end{figure}

%\input{examples}
%\input{numerics}
%!TEX root = main.tex 

\section{Numerical illustration}

In this section, we illustrate the theoretical results we have obtained previously.
In particular, we characterize below the shape of stability and unstability regions
for several different examples. 

We perform our numerical simulation in the setting of section \ref{subse motivation} 
using a trinomial tree (recombining) to approximate the Brownian motion and the terminal condition $\wh{\xi}= \cos(\wh{W}_T)$. 
Given a constant timestep $h>0$, the increment of the Brownian motion are
approximated by discrete random variables $\Delta \wh{W}_i$, $i \le n$, satisfying 
\begin{align*}
 \P(\Delta \wh{W}_i = \pm \sqrt{3{h}}) = \frac16 \quad \text{ and } \quad  \P(\Delta\wh{W}_i = 0) = \frac23\;.
\end{align*}
The $H$-coefficients are given by $H_i := \frac{\Delta \wh{W}_i}h$, $i \le n$, and are bounded.
Let us observe that, in the case of the implicit scheme $(\theta=1)$,
the stability condition of Proposition \ref{prop cond suff unidim} reads 
\begin{align}
 h \le \frac{1}{3 |L^Z|^2}\;. \label{eq num stab trinomial implicit}
\end{align}

On the graphs below, we plot the value $|Y_0| \wedge 10$ for different values of the parameter $h$ and various specifications of $f$.
Contrary to Section \ref{subse motivation},
for a fixed $h$, we chose to run the algorithm using $n=300$, which implicitly sets $T$ to be large. 
Doing so allows us to observe more clearly the various regions of stability and unstability 
for the specific choice of $f$, $\hat{\xi}$ and $h$.

% $T$ will be big to get a sharp characterisation of the stability/unstability.

\subsection{Linear specifications of $f$}

On Figures \ref{fig astab_mix_fz5} and \ref{fig astab_implinyfz_fz5} below, we plot $|Y_0|\wedge 10$ for $f(y,z) = ay + 5z$, for each $(a,h) \in [-3,0]\times(0,2]$.
This quantity corresponds here to the truncated absolute error between the scheme and the true solution, which is approximatively  equal to $0$ as $T$ is large.

\begin{figure}[htb!] 
\includegraphics[height=6cm, width=\textwidth]{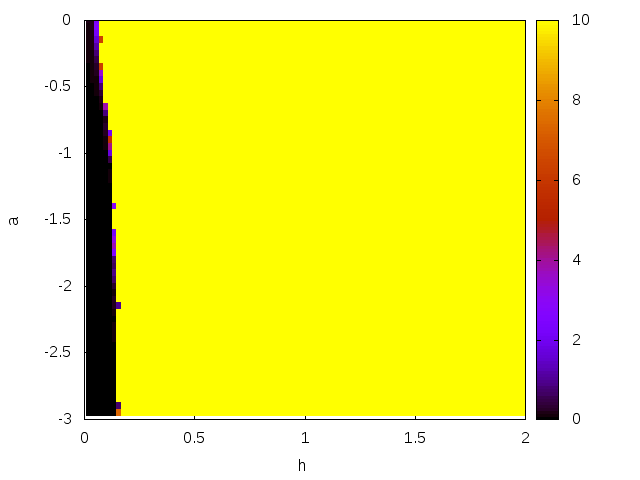} 
%[width = \textwidth]
\caption{Empirical stability of pseudo-explicit Euler scheme} 
\label{fig astab_mix_fz5}
\end{figure}

\vspace{2mm}
%We now compare the above result for the semi-explicit scheme with the one for the implicit scheme, which is clearly better.
\begin{figure}[htb!] 
\includegraphics[height = 6cm, width = 
\textwidth]{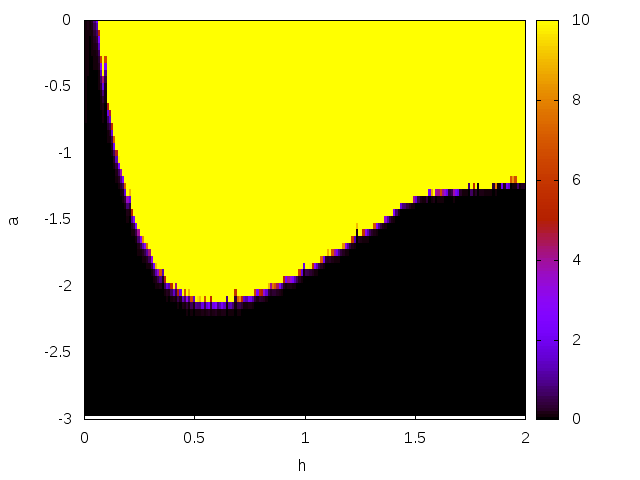} %[width = \textwidth]
\caption{Empirical stability of implicit Euler scheme} 
\label{fig astab_implinyfz_fz5}
\end{figure}

On both graphs, we are able to observe a stability region (in black) and an unstability region (in yellow). The shape of these regions is consistent
with the theoretical ones derived in Section \ref{se VN}, compare with Figures \ref{fig imp VN-stab region} and \ref{fig exp VN-stab region}. 

In Figure \ref{fig astab_imp_bZ}, we consider $f(z) =bz$ for $b \in [-5,5]$ and $h$ which varies between $0$ and $2$.

\begin{figure}[htb!] 
\includegraphics[height = 6cm, width = \textwidth]{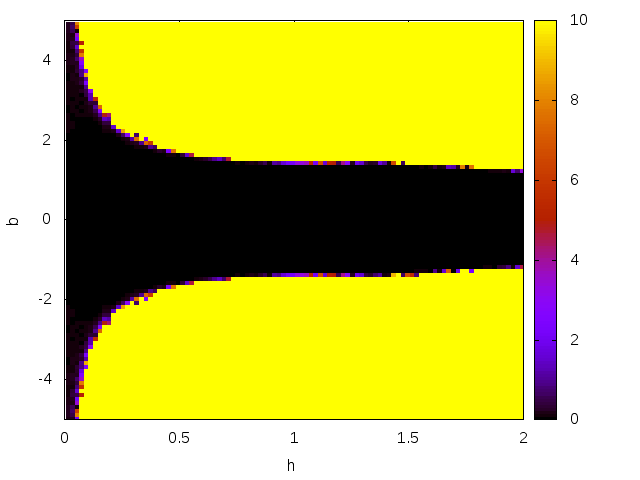} 
%[width = \textwidth]
\caption{Empirical stability of Euler scheme $f(y,z)=bz$} \label{fig astab_imp_bZ}
\end{figure}

\subsection{Non linear specifications of $f$}
In this section, we investigate the stability of the Euler scheme for some non-linear specification $z \mapsto f(z)$. 
%For a given $f$ and $\xi$, we do not know
%\emph{a priori} if the condition \ref{eq the condition} is necessary. The numerics below show that this is the case for some $f$ with our choice of $\xi$, see e.g. Figure \ref{fig astab_imp_babsZ}.

\subsubsection{$f(z)=b|z|$}
On the graph of Figure \ref{fig astab_imp_babsZ}, we plot the quantity $|Y_0| \wedge 10$ for $(b,h) \in [-5,5]\times(0,2]$.
In this example, we do not know the true value of $Y$ at time $t=0$. Nevertheless, we can clearly observe
the unstability region: in this case, the necessary condition seems to be related to \eqref{eq num stab trinomial implicit}.
%\textcolor{red}{tracer par dessus la courbe donnee par (4.1) ?}
% But we don't know the theoretical value, so we only plot $|Y|$ and truncate if $|Y|\ge10$...
% $Y$ does  not seem to be equal to $0$ in the stability zone...

\begin{figure}[htb!] 
\includegraphics[height = 6cm, width = 
\textwidth]{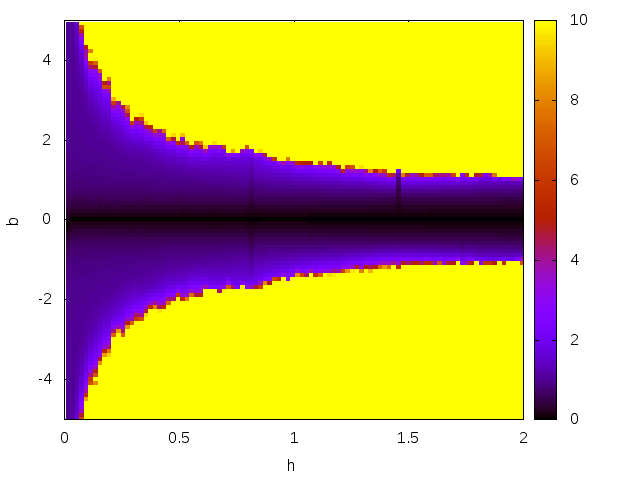} %[width = \textwidth]
\caption{Empirical stability of Euler scheme, $f(z)=b|z|$} \label{fig astab_imp_babsZ}
\end{figure}

\subsubsection{$f(z)=atan(bz)$}

On the plot in Figure \ref{fig astab_imp_atanbZ_HD}, we succesfully observe a stability region (in black) of the form predicted by \eqref{eq num stab trinomial implicit}. Outside this region, the behaviour of the algorithm seems fairly complicated. In particular, the algorithm is not robust outside the black (predicted) stability region as it seems to converge for some values of $h$ and not for others.

%\textcolor{magenta}{il faut choisir entre sin(bz) et atan(bz)}
\begin{figure}[htb!] 
\includegraphics[height = 6cm, width = 
\textwidth]{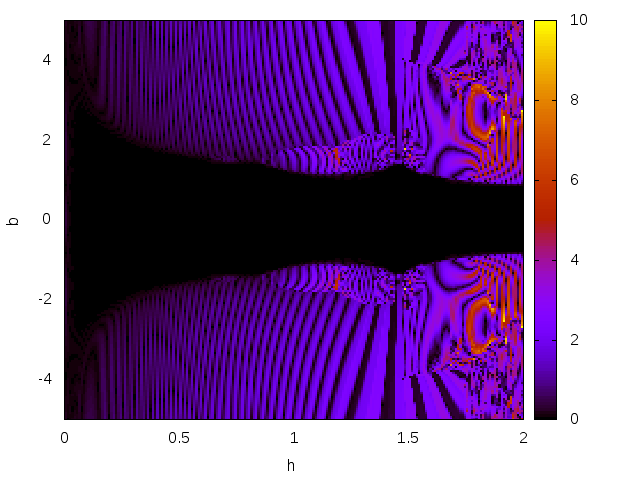} %[width = \textwidth]
\caption{Empirical stability of Euler scheme, $f(z)=atan(bz)$} \label{fig astab_imp_atanbZ_HD}
\end{figure}
\newpage

\def\cprime{$'$}

\end{document}